\documentclass[9pt,twocolumn]{IEEEtran}
\usepackage{amsmath,amssymb,euscript,graphicx,color,amstext,subcaption,parskip}
\usepackage{overpic}
\usepackage{tikz}
\usetikzlibrary{arrows}
\graphicspath{{./},{./figures/}}

\begin{document}
\newtheorem{thm}{Theorem}
\newtheorem{cor}[thm]{Corollary}
\newtheorem{conj}[thm]{Conjecture}
\newtheorem{lemma}[thm]{Lemma}
\newtheorem{prop}{Proposition}
\newtheorem{problem}[thm]{Problem}
\newtheorem{remark}[thm]{Remark}
\newtheorem{defn}[thm]{Definition}
\newtheorem{ex}[thm]{Example}

\newcommand{\mR}{{\mathbb R}}
\newcommand{\cR}{{\mathcal R}}
\newcommand{\cD}{{\mathcal D}}
\newcommand{\bC}{{\mathbb C}}
\newcommand{\bH}{{\mathbb H}}
\newcommand{\bS}{{\mathbb S}}
\newcommand{\derivt}[1]{\dfrac{\mathrm{d}#1}{\mathrm{dt}}}

\newcommand{\f}{{\mathfrak f}}

\newcommand{\diag}{\operatorname{diag}}
\newcommand{\trace}{\operatorname{trace}}
\newcommand{\spanof}{\operatorname{span}}
\newcommand{\argmin}{\operatorname{argmin}}
\newcommand{\ra}{\rangle}
\newcommand{\la}{\langle}
\newcommand{\be}{\begin{equation}}
\newcommand{\ee}{\end{equation}}

\newcommand{\ignore}[1]{}

\def\spacingset#1{\def\baselinestretch{#1}\small\normalsize}
\spacingset{.98}

\newcommand{\mike}{\color{magenta}}
\definecolor{grey}{rgb}{0.6,0.6,0.6}
\definecolor{lightgray}{rgb}{0.97,.99,0.99}

\title{Regularization and Interpolation\\
of Positive Matrices\thanks{This project was supported by the Air Force Office of Scientific Research (AFOSR) grants FA9550-12-1-0319 and FA9550-15-1-0045, NSF grant ECCS-1509387, NIH grants P41-RR-013218 and P41-EB-015902, and a postdoctoral fellowship through Memorial Sloan Kettering Cancer Center. 
}}

\author{Kaoru Yamamoto, Yongxin Chen, Lipeng Ning\\\hspace*{.2in}Tryphon T. Georgiou, and Allen Tannenbaum\thanks{K.\ Yamamoto is with the Department of Electrical and Computer Engineering, University of Minnesota, MN; email: kyamamot@umn.edu}
\thanks{Y.\ Chen is with the Department of Medical Physics, Memorial Sloan Kettering Cancer Center, NY; email: chen2468@umn.edu}
\thanks{L.\ Ning is with Brigham and Women's Hospital (Harvard Medical School), MA.; email: lning@bwh.harvard.edu}
\thanks{T.\ T. Georgiou is with the Department of Mechanical and Aerospace Engineering, University of California, Irvine, CA; email: tryphon@uci.edu}
\thanks{A.\ Tannenbaum is with the Departments of Computer Science and Applied Mathematics \& Statistics, Stony Brook University, NY; email: allen.tannenbaum@stonybrook.edu}}


\maketitle
{\begin{abstract}
We consider certain matricial analogues of optimal mass transport of positive definite matrices of equal trace. The framework is motivated by the need to devise a suitable geometry for interpolating positive definite matrices in ways that allow controlling the apparent tradeoff between ``aligning up their eigenstructure'' and ``scaling the corresponding eigenvalues''. Indeed, motivation for this work is provided by power spectral analysis of multivariate time series where, linear interpolation between matrix-valued power spectra generates push-pop artifacts. Push-pop of power distribuion is objectionable as it
corresponds to unrealistic response of scatterers.\\
\end{abstract}}


\section{Introduction}

The present paper is an attempt to develop a suitable matricial analogue of  
optimal mass transport (OMT). The basic problem of OMT refers to seeking a map or, in other words, a transportation plan that carries a given probability distribution to 
another in such a way that a certain specified transportation cost is minimized \cite{Rachev,Vil03,Vil08}. 
The original formulation of the problem by Monge in 1871 was motivated by civil engineering considerations, namely to transport dirt so as to level the ground. The problem achieved significant fame and notoriety due to inherent technical difficulties which persisted until the 1940's, at which time Kantorovich presented a 
relaxation of OMT in the form of a linear program. The relevance of this problem in the broader setting of resource allocation was already widely recognized and the impact secured a Nobel prize in Economics for Kantorovich in 1975. A new transformative phase of development in optimal mass transport began in the 1990's \cite{McC97,GanMcc96,JorKinOtt98,BenBre00} motivated by multitude of applications in physics, probability theory, image analysis, optimal control etc. In fact, OMT can be seen as an optimal control problem and stochastic formulations ensued (see \cite{Chen1,Chen2,Chen3} and the references therein).

In a recent publication \cite{ning2015matrix}, a matrix-valued version of OMT was put forth motivated by problems in spectral analysis of multivariable time series. More specifically, recall that the spectral content of scalar (slowly time-varying) time series is often displayed in the form of a spectrogram (time-frequency power distribution). While the spectrogram is illuminating in reflecting changes, it suffers from the perenial trade-offs between variability and resolution (uncertainty principle of Fourier methods) thereby necessitating spatio-temporal smoothing and regularization.
OMT is especially suited for this because the corresponding distance between density functions -- the Wasserstein metric, is ``weakly-continuous,'' in that small perturbations correspond to small changes in computed moments and vice versa. There has been no analogue for matrix-valued power spectra in the controls and signal processing literature. In fact,  \cite{ning2015matrix} was perhaps a first attempt which has its basis in Kantorovich's idea of seeking a joint density (matrix-valued in our case) in a suitable 
product space. In the present paper we explore alternative versions that are more deeply rooted in control ideas. To this end we seek dynamical equations that allow rotation of eigenvalues as well as scaling of the corresponding eigenvalues so as to generate paths between end-point matrices. The goal is to ensure that the optimal trajectory that connects end-points, promotes either rotation of eigenvectors or the scaling of eigenvalues, depending on the choice of regularization parameters.

For the development of a ``non-commutative'' version of OMT, i.e., a matrix-OMT, we replace probability 
density functions by density matrices $\rho$ -- the terminology ``density matrices'' is borrowed from quantum mechanics; these are positive semidefinite 
with unit trace. Now, ``transport'' corresponds to flow on the space of such 
matrices that minimizes a suitable cost functional. The insight and techniques gained
are aimed towards interpolating or regularizing sample covariances as well as matrix-valued power spectral  densities of multivariate time series -- they both reflect on how power varies with direction. 
In particular, we decompose the tangent space of the cones of positive definite matrices 
into two subspaces, one corresponding to rotating the 
eigenstructure and another to scaling the eigenvalues. The dynamics that will allow us to morph
matrices from a starting value to an end-point value are inspired by quantum mechanics.
Thus, in order to make the paper self-contained, we include some brief exposition of basic facts from quantum mechanics upon which we draw insight for the needed geometry. We also refer the interested reader to \cite{chen} for a parallel alternative formulation of matrix OMT, which once again, draws on the connection with the non-commutative geometry of quantum mechanics.

The paper is structured as follows. In Section~\ref{sec:quantum}, we present certain basic ideas 
of quantum mechanics that inspire the material in the paper. In Section~\ref{sec:tangentspace} we study
the tangent space of the cone of positive definite matrices. 
This leads to ideas on what are suitable cost functionals that promote a judicious balance between rotating the eigenstructure and scaling the eigenvalues (Section \ref{sec:interpolation}).
Numerical examples are given and discussed in Sections~\ref{sec:casestudy} and \ref{sec:regularization}.
Concluding remarks are provided in Section~\ref{sec:conclusion}.

\noindent
{\bf Notation:}
We denote by $\bH$ the set of $n \times n$ Hermitian matrices, 
$\bS$ the set of $n \times n$ skew-Hermitian matrices, and $\bH_+$ the cones of positive-semidefinite 
matrices. Since matrices are $n \times n$ throughout the paper, we will not explicitly 
note dependence on $n$. The commutator of two square matrices $A,B$ is denoted by 
$[A, B] := AB -BA$ and the anti-commutator by $\{A,B\} := AB+BA$.

\section{Quantum insights}\label{sec:quantum}

The development below draws concepts and insights from quantum mechanics and, therefore, we begin with a brief expository account of relevant basic facts. Detailed accounts can be sought in standard references, e.g., \cite{Sigal}.

\subsection{Schr\"odinger equation}
The evolution of closed quantum systems, i.e., one having no interaction with 
external quantum systems, is given by the 
\textit{time-dependent Schr\"odinger equation} \cite{Sigal}:
\be \label{eq:schroedinger}
\frac{\partial \psi}{\partial t} = X \psi
\ee
where $\psi \in {\mathbb C}^n$ and $X$ is a skew Hermitian matrix\footnote{More generally, $\psi$ belongs to a Hilbert space and accordingly $X$ is a skew Hermitian operator on that same Hilbert space. Typically $X$ is expressed as $-\frac{i}{\hbar}H$ where 
$H$ is a Hamiltonian (Hermitian) operator and $\hbar$ is the reduced Plank constant.}.  
Equation~\eqref{eq:schroedinger} describes a unitary evolution for the wave function, in our case, vector $\psi$; the quantum system is in a pure state. More generally, a system in a mixed state is described by the density matrix
\[
\rho= \sum_k \lambda_k\psi_k\psi_k^*
\]
with $\sum_k \lambda_k=1$ and evolves according to 
\be \label{eq:rotation}
\frac{\partial\rho}{\partial t}=X\rho-\rho X =: [X,\rho].
\ee
It is evident that if the system is in a pure state, it remains so, as the rank of
\begin{equation*}
\rho(t) = e^{Xt}\psi(0) \psi(0)^*e^{-Xt}
\end{equation*}
remains invariant. Likewise, if the system is in a mixed state
\begin{align*}
\rho(t)&= \sum_k \lambda_k(t)\psi_k(t)\psi_k^*(t)\\
&=e^{Xt}\left(\sum_k \lambda_k(t)\psi_k(0) \psi_k(0)^*\right)e^{-Xt},
\end{align*}
the eigenvalues $\lambda_k(t)$ of the density matrix remain invariant over time $t$, i.e.,
 $\lambda_k(t)=\lambda_k(0)\mbox{ for all }t.$
In other words, the evolution governed by \eqref{eq:rotation} only rotates in the same way the complete set of eigenvectors of the density matrix without changing the eigenvalues.

\subsection{Evolution of density matrices}

Decoherence and changes in the spectrum of $\rho$ are typically modeled through coupling with an \textbf{\it{ancilla}} which is another quantum system. The state of the original system is then obtained by {\it tracing out} the ancillatory component of the joint density operator. \textbf{\it{Lindblad's equation}} describes precisely such an evolution for the component of the original system. 
The Lindblad equation has the form
\be \nonumber \label{lindblad}
\frac{\partial \rho}{\partial t} = [X,\rho] - \sum_{k} \left(\frac12(Y_k\rho+\rho Y_k) - Z_k\rho Z_k^*\right)
\ee
where $Y_k=Z_k^* Z_k$. The presence of $-Z\rho Z^*$ ensures that $\trace(\rho)$ remains constant while both the eigenvalues and the eigenstructure may change over time.
Alternatively, one may consider more generally
\begin{equation}\label{eq:mf2}
\frac{\partial \rho}{\partial t}=[X,\rho] + u
\end{equation}
where $\trace(u)=0$ so as to preserve the trace of $\rho$. In fact,
in what follows, we will do exactly that and consider flows in directions corresponding to traceless component $u$. Positivity of the flow will be ensured as an added (convex) condition and will not be intrinsically encoded in $u$ (as in the Linblad equation where the right hand side is linear in $\rho$).

\section{Trace-preserving linear flow on positive matrices}\label{sec:tangentspace}
Consider the set of positive semidefinite matrices that are normalized to have trace one, 
\[
\cD :=\{\rho \in \bH_{+} \mid \trace(\rho)=1\}.
\]
As we noted earlier, we seek flows on $\cD$ that 
preserve trace.

The tangent space of $\bH_{+}$ at any point 
$\rho\in \bH_{+}$ is $\bH$. 
The subspace of traceless Hermitian components\footnote{Note that for $X\in\bS$ and $\rho\in\bH$, both $[X,\rho]^*=[X,\rho]$ and $\trace([X,\rho])=0$.}
\[
R_\rho:=\{[X,\rho] \mid X\in \bS\},
\]
is responsible for rotating the eigensubspaces of $\rho$ as we have noted in the previous section.

We now seek to identify
the orthogonal complement of $R_\rho$ so as to isolate the two directions that are responsible for rotation of eigenvectors and scaling of eigenvalues.
To this end, consider $u\in \bH$ such that
\begin{equation}\label{orthogonality}
\trace(u [X,\rho])=0 \quad \forall X\in \bS. 
\end{equation}
Since the trace is invariant under cyclic permutations,\begin{equation}\label{commutativity}
\trace([u, \rho]X)=0 \quad \forall X\in \bS.
\end{equation}
But $[u,\rho]$ is already in $\bS$, hence it is zero and therefore $u$ must commute with $\rho$.
Thus, from \eqref{orthogonality} we have that
 $u$ is in the orthogonal complement 
of $R_\rho$. We summarize our conclusion as follows.
\begin{prop} The tangent space $T_\rho$ of $\cD:=\{\rho \in \bH_{+} \mid \trace(\rho)=1\}$ 
at $\rho\in \cD$ 
can be decomposed as the direct sum
\[
T_\rho=R_\rho \oplus C_\rho
\]
of orthogonal components
\begin{align*}
R_\rho &= \left\{[X,\rho]\mid X\in \bS\right\} \mbox{ and}\\
C_\rho &= \left\{u \mid u\in \bH,~[u, \rho]=0,\mbox{ and }\trace(u)=0 \right\}.
\end{align*}
If $X(t)\in \bS$ and  $u(t)\in C_{\rho(t)}$ for all $t$, then $\trace(\rho(t))$ remains constant with $t$.
\end{prop}

\section{Interpolating flows between $\rho_0$ and $\rho_1$}\label{sec:interpolation}

Following on the previous rationale we may seek paths between density matrices 
$\rho_0$ and $\rho_1$, that minimize a suitable cost functional that allows trading off between 
the eigenstructure rotation specified by $X(t) \in \bS$ and the eigenvalue scaling affected by 
$u(t) \in \bH$. We are led to the following optimization problem:
\begin{equation*}
\begin{aligned}
{\text{minimize}}
& & \int_0^1 \left(\|X(t)\|_2+\epsilon \|u(t)\|_2\right)~\mathrm{dt} &\hspace*{15pt}\textbf{(Problem A)\quad} \\
 \text{subject to}
&& \dot \rho(t)= [X(t),\rho(t)]+u(t)\hspace*{20pt}\\
&&  \rho(0)=\rho_0, \;\rho(1)=\rho_1, \rho(\cdot)\geq 0\\
&& X\in \bS, \mbox{ and } u(t)\in C_{\rho(t)}.\hspace*{20pt}
\end{aligned}
\end{equation*}
Here $\epsilon$ is a choice of weight trading off the two alternative mechanisms for shifting eigenvalues and eigensubspaces to match the two end-point matrices. Also, here ``rotation'' and ``scaling'' may vary over time. 
Schematically the two operations are shown in the figure below (for a constant $X$).
{\small
\begin{center}
\begin{tikzpicture}[->,>=stealth',auto,node distance=6cm,
  thick,main node/.style={circle,draw,font=\sffamily\bfseries}]
 \node [main node] (rho0) {$\rho_0$};
 \node [main node](rho1) [right of=rho0]{$\rho_1$};
\draw [->, thick] (rho0) to [bend left=45] node {rotation: $e^{Xt}\rho_0 e^{-Xt}$} (rho1);
\draw [->, thick] (rho0) to [bend right=45,below] node {scaling: $\rho_0+\int^t_0u(t)~\mathrm{dt}$}(rho1);
\end{tikzpicture}
\end{center}
}	

Assuming a constant rate of rotation and a constant ``drift'' of the spectrum, we formulate a suitable problem as follows: 
\begin{equation} \label{u}
\begin{aligned}
{\text{minimize}}\hspace*{15pt}
 & \|X\|_2+\epsilon\|Z\|_2 \hspace*{65pt}\textbf{(Problem B)\quad}  \\
 \text{subject to}\hspace*{15pt}
& \dot \rho(t)= [X,\rho(t)]+\underbrace{e^{Xt}Ze^{-Xt}}_{u(t)}\\
 & \rho(0)=\rho_0, \;\rho(1)=\rho_1,\\
& X\in \bS, Z\in \bH, [\rho_0,Z]=0, \trace(Z)=0.
\end{aligned}
\end{equation}
Indeed, rotation and scaling of eigenvalues follow a constant ``drift'' and the solution
to
\[
\dot \rho(t)= [X,\rho(t)]+e^{Xt}Ze^{-Xt}
\]
is given by
\[
\rho(t)=e^{Xt}(\rho_0 + Zt)e^{-Xt}.
\]
We readily verify that $u(t)$ in \eqref{u} commutes with $\rho(t)$ as long as $Z$ commutes with $\rho_0$.

\section{Example: Interpolation of density matrices}\label{sec:casestudy}

In this section we highlight how interpolation is effected, via solving Problem B, as a proof of concept. Starting from two end-point density matrices, the framework allows constructing alternative paths connecting the two where one may tradeoff the two possible ways that the transition from one to the other may take place, i.e., allowing for the eigenvalues to adjust by ``scaling'' and the eigenvectors to ``rotate,'' respectively.

\begin{figure}[ht]
\centering
\begin{subfigure}[b]{0.45\textwidth}
\centering 
\includegraphics[scale=0.4]{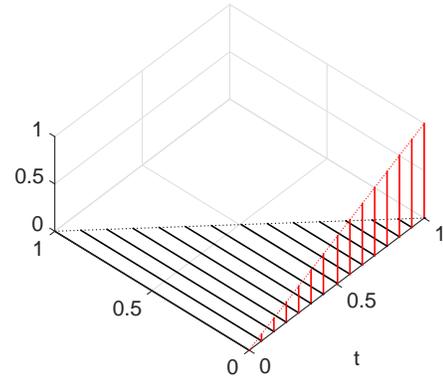}
\caption{Interpolation via adjusting eigenvalues: eigenvectors scaled according to corresponding eigenvalues.}
\label{fig:2d_scal}
\end{subfigure}
\begin{subfigure}[b]{0.45\textwidth}
\centering
\includegraphics[scale = 0.4]{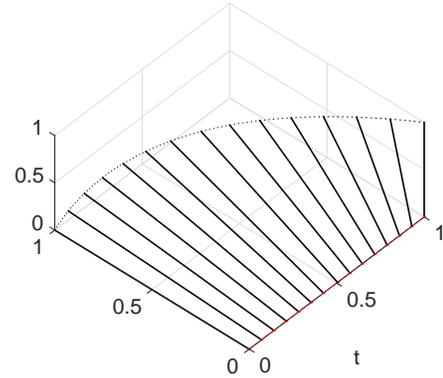}
\caption{Interpolation via rotating eigenvectors: eigenvectors scaled according to corresponding eigenvalues.}
\label{fig:2d_rot}
\end{subfigure}
\caption{Solutions are obtained by solving Problem B.}
\label{fig:2d_eigflow1}
\end{figure}

Consider the two density matrices
\[
\rho_0 = {\scriptsize \begin{bmatrix} 1&0\\0&0\end{bmatrix}}\mbox{ and }\rho_1 = {\scriptsize \begin{bmatrix} 0&0\\0&1\end{bmatrix}}.\] 
On one end, a choice of $\epsilon$ (vanishingly small) in Problem B leads to a path that displays a fade-in/fade-out effect of scaling the eigenvalues, so as to connect the two end-points (Fig.~\ref{fig:2d_scal}). No rotation of eigenvalues takes place. On the other end, for $\epsilon$ large, 
we obtain a path where ``rotation'' of the eigenvectors is less costly (Fig.~\ref{fig:2d_rot}); it is worth noting that in this case, since both matrices have rank one, the path remains rank one.


Motivation for our framework stems from multivariable time series analysis where power is often associated (e.g., in sensor arrays or radar) with the position of dominant scatterers. Fade-in-fade-out effects when interpolating or smoothing multivariable spectra are obviously undesirable as they create artifacts. Such fade-in-fade-out effects may be erroneously interpreted as due to the presence of additional scatterers beyond those that are present. The above rudimentary example may correspond to the case of two sensors reading a constant-frequency echo from a scatterer that changes its relative position with respect to the two. When recorded signals are correlated, the matrix-valued power spectrum at the corresponding frequency has (approximately) rank one. Likewise, movement of the scatterer that corresponds to a path between the two matrices, ought to have rank one. This exemplifies the need for paths that avoid push-pop for the corresponding eigenvalues (as linear interpolation would -- one eigenvalue reducing while another increasing at the same time).

\begin{figure}[th]
\centering
\begin{overpic}[trim = 20mm 0mm 20mm 0mm, clip, scale=0.7]{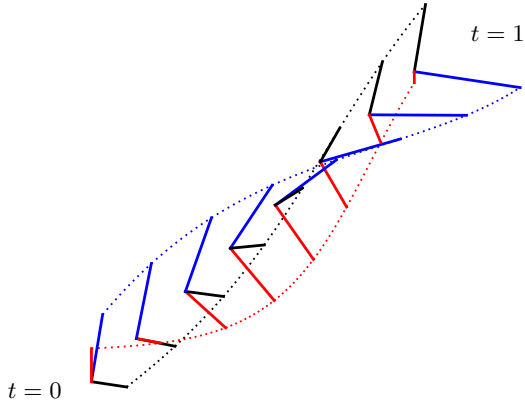}
\put(9,4){$t=0$}
\put(75,55){$t=1$}
\end{overpic}
\caption{
Solutions are obtained by solving Problem B.}

\label{fig:3d_eigflow}
\end{figure}

We close the section with an example of density matrices of higher dimension ($3$ in this case). Figure~\ref{fig:3d_eigflow} shows a path between two matrices 
\[
\rho_0 ={\scriptsize \begin{bmatrix} 1&0&0\\ 0&2&0 \\ 0&0&3\end{bmatrix}} \mbox{ and }{\rho_1 = \scriptsize\begin{bmatrix} 3&0&0\\ 0&2&0 \\ 0&0&1 \end{bmatrix}},
\]
for small $\epsilon$.
The solution is $\rho(t)=e^{Xt}(\rho_0 + Zt)e^{-Xt}$
 for
\[X = {\scriptsize \begin{bmatrix} 0 & 0 & 2.2\\ 0 & 0 & -2.2 \\-2.2& 2.2 &0 \end{bmatrix}} \mbox{ and }Z = {\scriptsize \begin{bmatrix} 1 & 0 & 0\\ 0 & 1&0\\ 0&0&-2 \end{bmatrix}}.\]
 This once again demonstrates the significance of penalizing rescaling of the eigenvalues which leads to adaptation of the corresponding eigenstructure via rotation. 

\section{Example: Regularization of noisy paths}\label{sec:regularization}
\begin{figure}[h]
\centering
\begin{subfigure}[b]{0.45\textwidth}
\centering 
\includegraphics[scale=0.4]{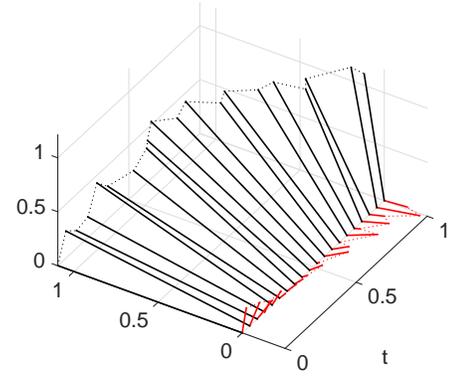}
\caption{Data: noisy matrices $\tilde{\rho}(t_i)$.}
\end{subfigure}
\begin{subfigure}[b]{0.45\textwidth}
\centering
\includegraphics[scale = 0.4]{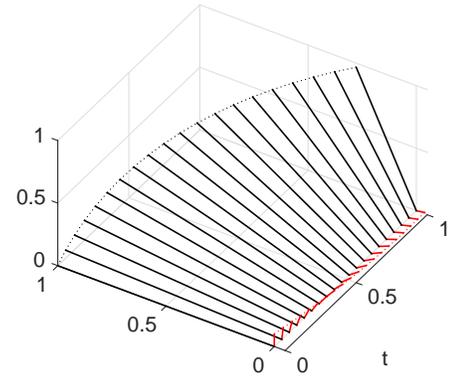}
\caption{Regularized path $\rho(t_i)$.}
\end{subfigure}
\caption[Text excluding the matrix]{Regularization of noisy matricial data: eigenvectors scaled according to corresponding eigenvalues.}
\label{fig:regularization}
\end{figure}

Besides interpolation problems, i.e., finding a path for $\rho(t)$ 
connecting two density matrices $\rho_0$ and $\rho_1$, the approach allows solving regularization problems where a smooth path is constructed to smooth out noisy measurements.
More specifically, 
given a noisy data set 
\[
\{\tilde{\rho}(t_i) \mid 0\leq t_1\leq \ldots \leq t_N \leq 1,\}
\]
we seek a 
smooth path $\rho(t)$ that approximately fits the data in a suitable sense. 
The key is to parametrize the path in a way consistent with the two ``orthogonal'' actions of rotating eigenvectors and scaling eigenvalues (as both may be needed), and penalize one more (typically, scaling). To this end, we solve the following 
extention of Problem B: 
\begin{equation*}
\begin{aligned}
& \underset{\rho_0,X,Z}{\text{minimize}}
& & \sum_{i=1}^N \left\|e^{Xt_i}(\rho_0 + Zt_i) e^{-Xt_i} - \tilde{\rho}(t_i)\right\|_2\\
&\text{subject to}
& 
&\hspace*{-5pt}{X\in \bS,\,Z\in \bH,\,\rho(\cdot)\geq 0,\,[\rho_0,Z]=0}\\
&&&\mbox{and }\trace(Z)=0.
\end{aligned}
\end{equation*}
The outcome is shown in Figure.~\ref{fig:regularization}.
For the purposes of illustration, the data set $\tilde{\rho}(t_i)$ is generated by adding a symmetric matrix-valued 
(uniform) noise $w(t)$ to a nominal flow 
$e^{Xt}\rho_0e^{-Xt}$  for $t_i\in \{0.05,\,0.1,\,0.15,\ldots,1\}$ where 
\[
\rho_0 = {\scriptsize \begin{bmatrix} 1.0 & 0\\ 0 & 0.1 \end{bmatrix}} \mbox{ and }X = {\scriptsize \begin{bmatrix} 0 & -1.6\\ 1.6 & 0 \end{bmatrix}}.
\]

\section{Conclusions}\label{sec:conclusion}
We developed an approach to constructing flows on density matrices. This allows interpolation and regularization of sample covariances and estimated power spectra of multivariable time series.
The general approach is control theoretic in that we select the flow (tangent direction) that minimizes a suitable cost functional. The choice of functional allows trading off the two basic mechanisms (rotating eigenvectors vs.\ scaling eigenvalues). Judicious balance between aligning up the 
eigenstructure and scaling the eigenvalues is necessitated by the fact that one of the two mechanisms alone may not suffice in generic situations. 

\bibliographystyle{plain}

\end{document}